\documentclass{article}

\def\a{\alpha}
\def\b{\beta}
\def\l{\lambda}
\def\L{\Lambda}
\def\f{\rightarrow}
\def\fb{\hookrightarrow}
\def\c{\widehat}
\def\q{\forall}

\def\G{\Gamma}
\def\D{\Delta}
  
\def\v{\vdash}

\def\et{\vedge}
\def\int{\bigcap}

\def\et{\wedge}

\def\sur{\overline}
\def\sous{\underline}

\def\F{\displaystyle\frac}

\def\N{\ifmmode{\rm I\mkern-3.1mu N\mkern0.5mu}\else{\rm I\kern-.16em N\hskip0.5pt\ }\fi\relax}

\newcommand{\cqfd}{\hbox{}\nobreak\hfill\nobreak$\spadesuit$}

\newtheorem{theoreme}{Th\'eor\`eme}[section]
\newtheorem{lemme}{Lemme}[section]
\newtheorem{corollaire}{Corollaire}[section]

\begin{document} 

\begin{center} 

{\Large\bf LES $I$-TYPES DU SYST\`EME ${\cal F}$}\\[0,5cm]

{\Large\bf $I$-TYPES OF SYSTEM ${\cal F}$}\\[0,3cm] 
 
{\large Karim NOUR\footnote{LAMA - Equipe de Logique -
Universit\'e de Chamb\'ery -
73376 Le Bourget du Lac -
Email nour@univ-savoie.fr}}
\end{center}

R\'esum\'e. - {\it Nous d\'emontrons dans ce papier que les types du
syst\`eme ${\cal F}$ habit\'es uniquement par des $\l I$-termes (les
$I$-types) sont \`a quantificateur positif. Nous pr\'esentons ensuite
des cons\'equences de ce r\'esultat et quelques exemples.}\\

Abstract. - {\it We prove in this paper that the types of system
${\cal F}$ inhabited uniquely by $\l I$-terms (the $I$-types) have a
positive quantifier. We give also consequences of this result and some
examples.}\\



\section{Introduction}

Le $\l I$-calcul est une restriction du $\l$-calcul o\`u on ne peut
abstraire sur une variable que si elle est libre dans le terme.  Dans
ce calcul, on a l'\'equivalence entre ``la normalisation faible'' et
``la normalisation forte''. H. Barendregt a d\'emontr\'e dans [1] que le $\l
I$- calcul est suffisant pour repr\'esenter les types de donn\'ees
courants et programmer toutes les fonctions calculables.\\

J.-Y. Girard a d\'emontr\'e dans [3] que l'on peut aussi repr\'esenter
les types de donn\'ees courants dans le syst\`eme de typage ${\cal
F}$. Cependant les fonctions repr\'esentables sur ces types sont
celles dont les preuves de terminaison se font dans l'arithm\'etique
de Peano du second ordre.\\

Les types du syst\`eme ${\cal F}$ qui repr\'esentent les types de
donn\'ees courants sont tous habit\'es par au moins un $\l$-terme qui
n'est pas un $\l I$-terme. Une question alors se pose : peut-on
repr\'esenter les types de donn\'ees courants par des types du
syst\`eme ${\cal F}$ habit\'es uniquement par des $\l I$-termes (ces
types sont appel\'es des {\it $I$-types})?\\

Le but de cet article est d'\'etudier les $I$-types du syst\`eme
${\cal F}$. On d\'emontre que tout quantificateur du second ordre d'un
$I$-type est positif. On donne ensuite une application de ce
r\'esultat sur les types entr\'ees et les types sorties d\'efinis dans
[2]. On montre aussi que pour v\'erifier si un type est un $I$-type,
on peut se liminter aux types simples (types sans quantificateurs). On
pr\'esente enfin quelques exemples de $I$-types.

\section{Notations et d\'efinitions}

\subsection{Quelques d\'efinitions du $\l$-calcul}

\'Etant donn\'es des $\l$-termes $t, u, u_1,..., u_n$, l'application de
$t$ \`a $u$ sera not\'ee $(t)u$, et $(... ((t)u_1) ...)u_n$ sera not\'e $(t)u_1...u_n$. \\

Si $t$ est un $\l$-terme, on d\'esigne par $Fv (t)$ l'ensemble de ses variables libres. \\

Si $u$ et $v$ sont des $\l$-termes, alors on note $\langle u,v \rangle$ le $\l$-terme $\l x
(x)uv$ o\`u $x$ ne figure pas dans $u$ et $v$.\\

On note ${\bf id} = \l x x$, ${\bf 0} = \l x \l y y$ et ${\bf 1}
= \l x \l y x$. L'{\it entier de Church} $\sous{n}$ est le $\l$-terme
$\l x \l f (f)...(f)x$ ($(f)$ r\'ep\'et\'e $n$ fois).\\  

On note $\f_{\b}$ (resp. $\f_{\b \eta}$) la $\b$-r\'eduction (resp. la $\b \eta$-r\'eduction).\\

Un $\l$-terme est dit {\it $\b$-normal} (resp. {\it $\b \eta$-normal}) s'il
ne contient pas de $\b$-redex (ni de $\b$-redex ni de $\eta$-redex).     \\

Un $\l$-terme est dit {\it r\'esoluble} si sa r\'eduction de t\^ete
termine.

\subsection{Le $\l I$-calcul}

{\bf D\'efinitions} 1) L'ensemble des {\it $\l I$-termes} (not\'e $\L I$) est d\'efinie par induction :

-- Si $x$ est une variable, alors $x \in \L I$.

-- Si $u,v \in \L I$, alors $(u)v \in \L I$.

-- Si $u \in \L I$ et $x$ est libre dans $u$, alors $\l x u \in \L I$.

2) Un {\it $\l K$-terme} est un $\l$-terme qui n'est pas un $\l I$-terme.\\

Le $\l I$-calcul poss\`ede les propri\'et\'es suivantes (voir [1]).

\begin{theoreme} 
1) Si $t \in \L I$, et $t \f_{\b \eta} t'$, alors $t' \in \L I$ et $Fv(t) = Fv(t')$.\\ 
2) Un $\l I$-terme est fortement normalisable ssi il est faiblement normalisable.
 \end{theoreme}

{\bf Notation} On note $\sur{0} = \l x \l f (((x){\bf id}){\bf id})f$
et, pour tout $n \in \N$, $\sur{n+1} = \sous{n+1}$.\\

Le th\'eor\`eme suivant montre qu'en $\l I$-calcul on peut
repr\'esenter toutes les fonctions partielles r\'ecursives (voir [1]).

\begin{theoreme} 
Pour toute fonction partielle r\'ecursive $f : \N^k \f \N$, il existe
un $\l I$-terme $\sur{f}$ tel que :  

 -- $(\sur{f})\sur{n_1}...\sur{n_k} \f_{\b} \sur{f(n_1,...,n_k)}$ si
$f(n_1,...,n_k)$ est d\'efinie. 

 -- $(\sur{f})\sur{n_1}...\sur{n_k}$ est non r\'esoluble sinon. 
 \end{theoreme}

\subsection{Le syst\`eme {$\cal F$}}

{\bf D\'efinitions} 1) Les {\it types} du syst\`eme {$\cal F$} sont
construits \`a partir des variables de type $X,Y,Z,...$ en utilisant
les op\'erations suivantes :

 -- Si $E$ et $F$ sont des types, alors $E \f F$ est un type.

 -- Si $E$ est un type, et $X$ est une variable de type, alors $\q X E$ est un type.

On d\'efinit d'une mani\`ere usuelle les variables {\it libres} et
 {\it li\'ees} d'un type. 

2) Soient $A,F$ deux types et $X$ une variable. Le type $A[F/X]$ est
   obtenu en remplacant dans $A$ toute occurrence de la variable $X$
   par le type $F$. 

3) Les r\`egles de typage du syst\`eme ${\cal F}$ sont les suivantes :  

\begin{center}
$(ax)$ $x_1 : A_1,...,x_n : A_n \v_{\cal F} x_i:A_i$ $(1\leq i\leq n)$
\\ [0.5cm]

$(\f_i)$ $\F{ \G,x:A \v_{\cal F} t:B } { \G\v_{\cal F} \l xt:A \f B }$ 
$\;\;\;\;\;\;\;\;\;\;\;\;\;\;\;\;\;\;\;\;\;\;\;\;\;$
$(\f_e)$ $\F{ \G\v_{\cal F} u:A \f B \quad \G\v_{\cal F} v:A} { \G\v_{\cal F}
(u)v:B }$ \\ [0.5cm]

$(\q_i)$ $\F{ \G\v_{\cal F} t:A \quad X{\rm \; non\; libre\; dans \;\G}} { \G\v_{\cal F} t:\q
XA }$
$\;\;\;\;\;\;$
$(\q_e)$ $\F{ \G\v_{\cal F} t:\q XA \quad G {\rm \;est\; un\; type}} { \G\v_{\cal F} t:A[G/X] }$
\end{center}

On \'ecrit $\G \v_{\cal F} t : A$ si $t$ est typable de type $A$ dans
le contexte $\G$.\\

Le syst\`eme $\cal F$ poss\`ede les propri\'et\'es suivantes (voir [3] et [4]).

\begin{theoreme} 
1) Un type est pr\'eserv\'e durant une $\b$-r\'eduction.

2) Un $\l$-terme typable est fortement normalisable.
 \end{theoreme} 

{\bf D\'efinition} Un type $D$ du syst\`eme ${\cal F}$ est dit {\it
propre} ssi si $\q X E$ est un sous-type de $D$ alors $X$ est libre
dans $E$.\\

 Dans la suite on restreint le syst\`eme ${\cal F}$ aux types
 propores. Pour cela on modifie l\'eg\`erement la r\`egle $(\q_e)$ en
 demandant \`a $X$ d'\^etre libre dans $A$.\\

Les lemmes suivants seront utilis\'es dans la suite (voir [3] et [4]).\\

{\bf D\'efinition} Soient $Id = \q X \{ X \f X \}$, $Bool = \q X \{ X
\f (X \f X) \}$ et $Ent = \q X \{ X \f [(X \f X) \f X]\}$.

\begin{lemme} Soit $t$ un $\l$-terme $\b$-normal.

1) $\v_{\cal F} t : Id$ ssi $t = {\bf id}$.

2) $\v_{\cal F} t : Bool$ ssi $t = {\bf 0}$ ou $t = {\bf 1}$.

3) $\v_{\cal F} t : Ent$ ssi il existe $n \in\N$ tel que $t = \sous{n}$ ($n \geq 0$).
\end{lemme}

{\bf D\'efinition} Si $A$ et $B$ sont deux types alors on note $A \et B$ le type $\q X \{(A \f (B
\f X)) \f X \}$ o\`u $X$ est une variable qui ne figure pas dans $A$ et $B$.

\begin{lemme} 
1) Si $\G \v_{\cal F} u : A$ et $\G \v_{\cal F} v: B$, alors $\G
   \v_{\cal F} \langle u,v\rangle : A \et B$.

2) Si $\G \v_{\cal F} t : A \et B$, alors $\G \v_{\cal F} (t){\bf 1} : A$ et 
$\G \v_{\cal F} (t){\bf 0} : B$.
\end{lemme} 
  
{\bf D\'efinition} Un {\it contexte} $C \langle \; \rangle$ est un
terme du $\l$-calcul avec un trou (une seule occurrence d'une constante
sp\'eciale du $\l$-calcul). On
note  $C \langle u \rangle$ le r\'esultat de la substitution du trou
de $C \langle \; \rangle$ par $u$ et ceci sans renommage des variables
li\'ees de $C \langle \; \rangle$. 

\begin{lemme} Si dans le typage $\G \v_{\cal F} C \langle u \rangle : A$ on a 
$\G' \v_{\cal F} u : B$, alors pour tout $\l$-terme $v$ tel que $\G' \v_{\cal F} v : B$ on a
$\G \v_{\cal F} C \langle v \rangle : A$
\end{lemme} 

{\bf Preuve} Il suffit de remplacer dans le typage $\G \v_{\cal F} C
\langle u \rangle : A$ l'arbre de typage de $\G' \v_{\cal F} u : B$
par celui de $\G' \v_{\cal F} v : B$. \cqfd \\

{\bf Notations} Pour simplifier on note la formule $\q X_1 ... \q X_n
F$ par $\q \sur{X} F$ et la formule $A_1 \f (A_2 \f ...(A_n \f A)...)$
par $A_1,A_2,...,A_n \f A$.

\section{Les $I$-types du syst\`eme ${\cal F}$}

{\bf D\'efinition} Un type clos $D$ du syst\`eme ${\cal F}$ est dit {\it $I$-type} ssi si $t$ est un
$\l$-terme clos $\b$-normal tel que $\v_{\cal F} t : D$, alors $t$ est un $\l I$-terme.\\

Nous allons d\'emontrer qu'un $I$-type du syst\`eme ${\cal F}$
(habit\'e par au moins un $\l$-terme) est un type \`a quatificateur
positif. Nous pr\'esentons tout d'abord, sur des exemples, la m\'ethode adopt\'ee.\\

{\bf Exemples} Soient $E = \q X \{ \q Y (Id \f Y) \f Id \}$ et $F = \q X \{ \q Y (Y \f Id) \f Id
\}$. Ces types ne sont pas \`a quantificateur positif car le quantificateur $\q Y$ occure
n\'egativement dans $E$ et $F$. Nous allons montrer comment fabriquer \`a partir d'un $\l
I$-terme de type $E$ (resp. de type $F$) un $\l K$-terme de m\^eme type. Consid\'erons les deux typages suivants : 
\begin{center}
$\F{\F{\F{
\F{x : \q Y (Id \f Y) \v x : \q Y (Id \f Y)}
{x : \q Y (Id \f Y) \v x : Id \f Id} 
\quad
\F{}{\v_{\cal F}{\bf id} : Id}}
{x : \q Y (Id \f Y) \v (x){\bf id} : Id}}
{\v \l x (x){\bf id} : \q Y (Id \f Y) \f Id}}
{\v \l x (x){\bf id} : E}$
\end{center}

\begin{center}
$\F{\F{\F{
\F{x : \q Y (Y \f Id) \v x : \q Y (Y \f Id)}
{x : \q Y (Y \f Id) \v x : Id \f Id}
\quad
\F{}{\v_{\cal F}{\bf id} : Id}}
{x : \q Y (Y \f Id) \v (x){\bf id} : Id}}
{\v \l x (x){\bf id} : \q Y (Y \f Id) \f Id}}
{\v \l x (x){\bf id} : F}$
\end{center}

Les deux typages suivants donnent des $\l K$-termes de type $E$ et $F$ :

\begin{center}
$\F{\F{\F{\F{\F{x : \q Y (Id \f Y) \v x : \q Y (Id \f Y)}
{x : \q Y (Id \f Y) \v x : Id \f (Bool \f Id)}
\quad \F{}{\v {\bf id} : Id}
}
{x : \q Y (Id \f Y) \v (x){\bf id} : Bool \f Id}
\quad \F{}{\v {\bf 0} : Bool}
}
{x : \q Y (Id \f Y) \v (x){\bf id}\;{\bf 0} : Id}}
{\v \l x (x){\bf id}\;{\bf 0} : \q Y (Id \f Y) \f Id}}
{\v \l x (x){\bf id}\;{\bf 0} : E}$
\end{center}

\begin{center}
$\F{\F{\F{
\F{x : \q Y (Y \f Id) \v x : \q Y (Y \f Id)}
{x : \q Y (Y \f Id) \v x : (Bool \et Id) \f Id}
\quad
\F{\F{}{\v {\bf 0} : Bool} \quad \F{}{\v {\bf id} : Id}}{\v_{\cal F}
\langle {\bf 0},{\bf id}\rangle : Bool \et Id}}
{x : \q Y (Y \f Id) \v (x)\langle {\bf 0},{\bf id} \rangle : Id}}
{\v \l x (x)\langle {\bf 0},{\bf id} \rangle : \q Y (Y \f Id) \f Id}}
{\v \l x (x)\langle {\bf 0},{\bf id} \rangle : F}$
\end{center}

Remarquons que dans le premier typage on a remplac\'e la variable $Y$ (qui est en position
n\'egative dans $E$) par $Bool \f Id$ et dans le deuxi\`eme typage on a remplac\'e $Y$ (qui est en
position positive dans $F$) par $Bool \et Id$. On peut v\'erifier facilement qu'en remplacant $Y$
par $Bool \f (Bool \et Id)$ on trouve aussi  des $\l K$-termes de type $E$ et $F$. Signalons enfin
qu'il suffit de remplacer le type $Bool$ par un type quelconque
habit\'e par un $\l K$-terme ou par une nouvelle variable de type. Nous allons montrer qu'on peut g\'en\'eraliser cette m\'ethode pour un type
quelconque. Ceci n\'ecessite l'introduction de plusieurs notions et la d\'emonstartion de
plusieurs r\'esultats.\\

{\bf D\'efinition} Pour tout type $A$ et pour toute variable $X$, on
d\'efinit deux constantes du $\l$-calcul ${\cal U}_{A,X}$ et ${\cal
V}_{A,X}$. Le {\it $\l_{{\cal U}{\cal V}}$-calcul} est obtenu en
consid\'erant les r\`egles de r\'eduction suivantes :
\begin{center}
$(\l x u) v \fb_{\b} u[v/x]$
$\matrix{({\cal U}_{Y,X})t \fb_{u} t \; si \; Y \not = X & ({\cal
V}_{Y,X})t \fb_{v} t \; si \; Y \not = X \cr ({\cal U}_{B \f C,X})t
\fb_{u}  \l y ({\cal U}_{C,X})(t)({\cal V}_{B,X})y & ({\cal V}_{B \f
C,X})t \fb_{v} \l y ({\cal V}_{C,X})(t)({\cal U}_{B,X})y \cr ({\cal
U}_{\q Y B,X})t \fb_{u} ({\cal U}_{B,X})t & ({\cal V}_{\q Y B,X})t
\fb_{v} ({\cal V}_{B,X})t}$
\end{center}
On \'ecrit $t \fb t'$ si $t'$ est obtenu \`a partir de $t$ en appliquant un nombre fini de fois
les r\`egles pr\'ec\'edentes.

\begin{lemme}
La cl\^oture r\'eflexive et transitive de la r\'eduction $\fb_{u}$
(resp. de $\fb_{v}$) est fortement normalisable.
\end{lemme}

{\bf Preuve} On d\'efinit, par induction sur les $\l_{{\cal U}{\cal
V}}$-termes, une notion de longueur : $N(x) = 0$, $N((u)v) = N(u) +
N(v)$, $N(\l x u) = N(u)$ et $N({\cal U}_{A,X}) = N({\cal V}_{A,X}) =
L(A)$ o\`u $L(A)$ est le nombre des connecteurs logiques de $A$. Il
est clair que si $t \fb_{u} t'$ ou $t \fb_{v} t'$, alors $N(t) >
N(t')$. D'o\`u le r\'esultat.\cqfd \\

{\bf D\'efinitions} Soit $t$ un $\l_{{\cal U}{\cal V}}$-terme et $E$
un ensemble de variables.

1) Une occurrence d'une variable $x$ dans $t$ est dite {\it
E-inactive} ssi ou bien $x \in E$ ou bien il existe un sous-terme
$u$ de $t$ tel que $x \in Fv(u)$ et il existe une occurrence d'une
variable $E$-inactive $y$ telle que $((...(y) u_1 ...u_i)\l y_1 ...\l
y_m \l x u)u_{i+1}...u_n$ est un sous-terme de $t$.

2) Un sous-terme $u$ de $t$ est dit {\it $E$-inactif} ssi $u = (x)
u_1...u_n$ et $x$ est $E$-inactive dans $t$. 
 
3) Un sous-terme $u$ de $t$ est dit {\it $E$-passif} ssi il existe un $\l_{{\cal
U}{\cal V}}$-terme $E$-inactif $(x) u_1 ...u_n$ de $t$ tel que $u_i = \l y_1 ...\l y_m u$. 

\begin{lemme} Soient $u,v$ des $\l_{{\cal U}{\cal V}}$-termes, $E$ un ensemble de variables et
$x \not \in E$. Les occurrences des variables $E$-inactives (resp. les
termes $E$-inactifs, les termes $E$-passifs) de $u$ et $v$ sont des
occurrences des variables $E$-inactives (resp. des termes
$E$-inactifs, des termes $E$-passifs) de $u[v/x]$.
\end{lemme}

{\bf Preuve} Facile. \cqfd \\

{\bf D\'efinition} Soit $E$ un ensemble de variables. Un $\l_{{\cal
U}{\cal V}}$-terme $t$ est dit {\it $E$-bon} ssi chaque occurrence de
${\cal U}_{B,X}$ dans $t$ est appliqu\'ee \`a un seul argument $w$
et$({\cal U}_{B,X})w$ est $E$-passif et chaque occurrence de ${\cal
V}_{B,X}$ dans $t$ est appliqu\'ee \`a un $\l_{{\cal U}{\cal
V}}$-terme $E$-inactive.\\

Il est clair qu'un $\l_{{\cal U}{\cal V}}$-terme $E$-bon a l'une des
formes suivantes : $\l x u$, $(x)u_1 ... u_n$, $(\l x u)v v_1...v_n$
ou $({\cal V}_{B,X})(x)u_1 ... u_n$. Le lemme suivant caract\'erise
les $\l_{{\cal U}{\cal V}}$-termes $E$-bons.

\begin{lemme} 
1) $\l x u$ est $E$-bon ssi $u$ est $E-\{x\}$-bon.

2) $(x) u_1 ... u_n$ est $E$-bon ssi pour tout $(1 \leq i \leq n)$,
$u_i$ est $E$-bon ou (si $x \in E$) $u_i =\l y_1 ...\l y_m ({\cal
U}_{B,X})w$ et $w$ est $E-\{y_1,...,y_m\}$-bon.

3) $(\l x u) v v_1 ... v_n$ est $E$-bon ssi $v$ est $E$-bon, pour tout
$(1 \leq i \leq n)$, $v_i$ est $E$-bon et $u$ est $E-\{x\}$-bon.

4) $({\cal V}_{B,X})(x)u_1 ... u_n$ est $E$-bon ssi $x \in E$ et $(x) u_1 ... u_n$ est $E$-bon. 
\end{lemme}

{\bf Preuve} Facile. \cqfd

\begin{lemme} Soient $u,v$ des $\l_{{\cal U}{\cal V}}$-termes et $E$ un
ensemble de variables. Si $u,v$ sont $E$-bons, alors $(u)v$ est $E$-bon. 
\end{lemme}

{\bf Preuve} Facile.  \cqfd

\begin{lemme} Soient $u,v$ des $\l_{{\cal U}{\cal V}}$-termes, $E$ un ensemble de variables et
$x \not \in E$. Si $u,v$ sont $E$-bons, alors $u[v/x]$ est $E$-bon.
\end{lemme}

{\bf Preuve} On utilise le lemme 3.2.   \cqfd

 \begin{theoreme}
Si $t$ est un $\l_{{\cal U}{\cal V}}$-terme $E$-bon et $t \fb t'$,
alors $t'$ est $E$-bon.  
\end{theoreme}

{\bf Preuve} Il suffit de faire la preuve pour un seul pas de
r\'eduction. On proc\`ede par induction sur $t$ et on utilise les lemmes
3.3, 3.4 et 3.5.

 -- Si $t = \l x u$, alors il suffit d'appliquer l'hypoth\`ese
d'induction et utiliser le lemme 3.5.

 -- Si $t = (x) u_1 ... u_n$, alors on a deux cas a voir.  Si $u_i$
est $E$-bon et se r\'eduit a $u'_i$, alors il suffit d'appliquer
l'hypoth\`ese d'induction et utiliser le lemme 3.3.  Si $x \in E$,
$u_i =\l y_1 ...\l y_m ({\cal U}_{B,X})w$ et $w$ est
$E-\{y_1,...,y_m\}$-bon, alors le r\'esultat est \'evident si on fait
la r\'eduction dans $w$. Si on r\'eduit le ${\cal U}$-redex, alors
ceci d\'epend de $B$.

- Si $B = Y \not = X$, alors $u_i$ se r\'eduit \` a $\l y_1 ...\l y_m
w$. Or comme $w$ est $E-\{y_1,...,y_m\}$-bon, alors, en utilisant le
lemme 3.3, $\l y_1 ...\l y_m w$ est $E$-bon et donc $t'$ aussi.

- Si $B = C \f D$, alors $u_i$ se r\'eduit \`a $\l y_1 ...\l y_m \l
y({\cal U}_{D,X})(w)({\cal V}_{C,X})y$. Or comme $w$ est
$E-\{y_1,...,y_m\}$-bon, alors (par d\'efinition) $t'$ est $E$-bon.

- Si $B = \q Y C$, alors le r\'esultat est \'evident.

-- Si $t = (\l x u)v v_1...v_n$, alors le r\'esultat est \'evident si
on fait la r\'eduction dans les termes $u,v,v_1,...,v_n$.  Si $t' =
(u[v/x])v_1...v_n$, alors le r\'esultat provient des lemmes 3.4 et
3.5.

-- Si $t = ({\cal V}_{B,X})(x)u_1 ... u_n$, alors le r\'esultat est
\'evident si on fait la r\'eduction dans les termes $u_1,...,u_n$.  Si
on r\'eduit le ${\cal V}$-redex, alors ceci d\'epend de $B$.

- Si $B = Y \not = X$, alors $t' = (x)u_1 ... u_n$ qui est $E$-bon.

- Si $B = C \f D$, alors $t$ se r\'eduit \`a $\l y({\cal V}_{D,X})
(x)u_1 ... u_n ( {\cal U})y$. Donc (par d\'efinition) $t'$ est
$E$-bon.

- Si $B = \q Y C$, alors le r\'esultat est \'evident.\cqfd \\

{\bf D\'efinition} Un $\l_{{\cal U}{\cal V}}$-terme $t$ est dit {\it
bon} ssi il est $Fv(t)$-bon.

\begin{corollaire}
Si $t$ est un $\l_{{\cal U}{\cal V}}$-terme bon et $t \fb t'$, alors
$t'$ est bon.
\end{corollaire}

{\bf Preuve} D'apr\`es le th\'eor\`eme 3.1.\cqfd \\

{\bf D\'efinition} Soit ${\cal U}$ et ${\cal V}$ deux constantes
fix\'ees du $\l$-calcul. Pour tout type $A$ et pour toute variable
$X$, on d\'efinit par induction deux $\l$-termes ${\cal I}'_{A,X}$ et
${\cal J}'_{A,X}$ de la mani\`ere suivante : 
\begin{center}
$\matrix{{\cal I}'_{Y,X} = {\bf id} \;  si \; Y \not = X & 
{\cal J}'_{Y,X} = {\bf id} \; si \; Y \not = X \cr
{\cal I}'_{X,X} = {\cal U} & 
{\cal J}'_{X,X} = {\cal V} \cr
{\cal I}'_{B \f C,X} = \l x \l y ({\cal I}'_{C,X})(x)({\cal
J}'_{B,X})y &
{\cal J}'_{B \f C,X} = \l x \l y ({\cal J}'_{C,X})(x)({\cal I}'_{B,X})y \cr
{\cal I}'_{\q Y B,X} = \l x ({\cal I}'_{B,X}) x 
& {\cal J}'_{\q Y B,X} = \l x ({\cal J}'_{B,X}) x}$
\end{center}

\begin{lemme}
Pour tout type $A$ et pour toute variable $X$, les $\l$-termes ${\cal
I}'_{A,X}$ et ${\cal J}'_{A,X}$ sont des $\l I$-termes.
\end{lemme}

{\bf Preuve} Par induction sur $A$.  \cqfd \\

{\bf D\'efinition} On associe \`a chaque $\l_{{\cal U}{\cal V}}$-terme
$t$ un $\l$-terme not\'e $\c{t}$ de la mani\`ere suivante : $\c{x} =
x$, $\c{\l x u} = \l x \c{u}$, $\c{(u)v} = (\c{u})\c{v}$, $\c{{\cal
U}_{A,X}} = {\cal I}'_{A,X}$ et $\c{{\cal V}_{A,X}} = {\cal
J}'_{A,X}$.

\begin{lemme}
Si $u,v$ sont des $\l_{{\cal U}{\cal V}}$-termes, alors $\c{u[v/x]} =
\c{u}[\c{v}/x]$.
\end{lemme}

{\bf Preuve} Par induction sur $u$. \cqfd

\begin{lemme}
Soient $u,v$ des $\l_{{\cal U}{\cal V}}$-termes. Si $u \fb v$, alors
$\c{u} \f_{\b} \c{v}$.
\end{lemme}

{\bf Preuve} Par induction sur $u$ et on utilise le lemme 3.7.  \cqfd

\begin{lemme}
Si $t$ est un $\l_{{\cal U}{\cal V}}$-terme tel que $\c{t}$ est
fortement normalisable, alors $t$ est fortement
normalisable. \end{lemme}

{\bf Preuve} Sinon, alors il existe une suite $(t_i)_i$ telle que $t =
t_0$ et ($t_i \fb_{\b} t_{i+1}$ ou $t_i \fb_{u} t_{i+1}$ ou $t_i
\fb_{u} t_{i+1}$). D'apr\`es les lemmes 3.1 et 3.9, il existe une
suite croissante d'entiers $(n_i)_i$ telle que $n_0 = 0$ et
$\c{t_{n_i}} \f_{\b} \c{t_{n_{i+1}}}$. Ce qui contredit le fait que
$\c{t}$ est fortement normalisable. \cqfd \\

Avec les hypoth\`eses du lemme pr\'ec\'edent (et en utilisant le lemme
3.8) on a unicit\'e de la forme normale de $t$.\\

Donc la suite on fixe une constante du $\l$-calcul $\a$ et une constante de type $O$.\\

{\bf D\'efinitions} 1) Soient $U = \l x \l d \langle x , \a \rangle$ et $V = \l x
(x) \a {\bf 1}$. Pour tout type $A$ et pour toute variable $X$, on
note ${\cal I}_{A,X} = {\cal I}'_{A,X}[ U / {\cal U} , V / {\cal V}]$
et ${\cal J}_{A,X} = {\cal J}'_{A,X}[ U / {\cal U} , V / {\cal V}]$.

2) Pour tout type $G$ du syst\`eme ${\cal F}$, on note $G^{\circ} = O
\f (G \et O)$.

\begin{lemme}
Pour tout type $A$ et pour toute variable $X$, on a :

1) $\a : O \v_{\cal F} {\cal I}_{A,X} : \q Y \{ A[Y/X] \f
   A[Y^{\circ}/X] \}$.

2) $\a : O \v_{\cal F} {\cal J}_{A,X} : \q Y \{ A[Y^{\circ}/X] \f A[Y/X] \}$.
\end{lemme}

{\bf Preuve} Par induction sur $A$. 

-- Si $A = Y \not = X$, alors $A[Y/X] = A[Y^{\circ}/X] = Y$, ${\cal
I}_{A,X} = {\cal J}_{A,X} = {\bf id}$ et on a $\a : O \v_{\cal F} {\bf id}
: \q Y \{ Y \f Y \}$.

-- Si $A = X$, alors $A[Y/X] = Y$, $A[Y^{\circ}/X] = Y^{\circ}$,
${\cal I}_{A,X} = U$ et ${\cal J}_{A,X} = V$. On a les deux typages
suivants :

\begin{center}
$\F{\a : O , x : Y  , d : Y \v_{\cal F} \langle x,\a \rangle : Y \et O}
{\a : O \v_{\cal F} U : \q Y \{ Y \f Y^{\circ} \}}$

$\F{\a : O , x : Y^{\circ}  , \v_{\cal F} (x)\a{\bf 1} : Y }
{\a : O \v_{\cal F} V : \q Y \{ Y^{\circ} \f Y \}}$
\end{center}

-- Si $A = B \f C$, alors $A[Y/X] = B[Y/X] \f C[Y/X]$, et
$A[Y^{\circ}/X] = B[Y^{\circ}/X] \f C[Y^{\circ}/X]$. Par hypoth\`ese
d'induction on a $\a : O \v_{\cal F} {\cal I}_{C,X} : \q Y \{ C[Y/X]
\f C[Y^{\circ}/X] \}$ et $\a : O \v_{\cal F} {\cal J}_{B,X} : \q Y \{
B[Y^{\circ}/X] \f B[Y/X] \}$. Donc
\begin{center}
$\F{\F{\F{
\a : O , x : B[Y/X] \f C[Y/X] , y : B[Y^{\circ}/X]  \v_{\cal F} ({\cal J}_{B,X})y : B[Y/X]}
{\a : O , x : B[Y/X] \f C[Y/X] , y : B[Y^{\circ}/X]  \v_{\cal F} (x)({\cal J}_{B,X})y : C[Y/X]}}
{\a : O , x : B[Y/X] \f C[Y/X] , y : B[Y^{\circ}/X]  \v_{\cal F} ({\cal I}_{C,X})({\cal J}_{B,X})y : C[Y/X]}}
{\a : O , \v_{\cal F} {\cal I}_{A,X} : \q Y \{ A[Y/X] \f A[Y^{\circ}/X] \}}$
\end{center} 
De m\^eme on d\'emontre que $\v_{\cal F} {\cal J}_{A,X} : \q Y \{ A[Y^{\circ}/X] \f A[Y/X] \}$.

-- Si $A = \q Z B$, alors $A[Y/X] = \q Z B[Y/X]$, et $A[Y^{\circ}/X] =
\q Z B[Y^{\circ}/X]$. Par hypoth\`ese d'induction on a $\a : O
\v_{\cal F} {\cal I}_{B,X} : \q Y \{ B[Y/X] \f C[B^{\circ}/X]
\}$. Donc
\begin{center}
$\F{\F{
\a : O , x : \q Z B[Y/X]  \v_{\cal F} ({\cal I}_{B,X})x : B[Y^{\circ}/X]}
{\a : O , x : \q Z B[Y/X]  \v_{\cal F} ({\cal I}_{B,X})x : \q Z B[Y^{\circ}/X]}}
{\v_{\cal F} {\cal I}_{A,X} : \q Y \{ A[Y/X] \f A[Y^{\circ}/X] \}}$
\end{center} 

De m\^eme on d\'emontre que $\v_{\cal F} {\cal J}_{A,X} : \q Y \{
A[Y^{\circ}/X] \f A[Y/X] \}$. \cqfd \\ 

{\bf D\'efinitions} 1) Soit $t$ un terme $\b$-normal ; $t$ s'\'ecrit
$\l x_1 ... \l x_n (x) t_1 ...  t_m$. Les sous-termes {\it essentiels}
de $t$ sont, par d\'efinition, $(x) t_1 ... t_m$ et les sous-termes
essentiels des $t_i$ $( 1 \leq i \leq m)$.

2) Soient $\G$ un contexte, $A$ un type, et $t$ un $\l$-terme
$\b$-normal. On dit que $t$ est un $\l$-terme {\it $\eta$-long} de
type $A$ dans le contexte $\G$ ssi $\G \v_{\cal F} t : A$ et dans ce
typage tous les sous termes essentiels de $t$ sont typ\'es par des
variables de type.

\begin{lemme}
Si $t$ un $\l$-terme $\b$-normal tel que $\G \v_{\cal F} t : A$, alors
il existe un $\l$-terme $\eta$-long $t'$ de type $A$ dans le contexte
$\G$ tel que $t' \f_{\eta} t$.
\end{lemme}

{\bf Preuve} Il suffit de remplacer chaque sous-terme essentiel $u$ de
$t$ de type $B = \q \sur{X_0} (B_1 \f ...\q \sur{X_{n-1}} (B_n \f \q
\sur{X_n}X)...)$ par le $\l$-terme $\l x_1 ... \l x_n (u) x_n
... x_1$. \cqfd \\

{\bf D\'efinitions} 1) Un sous-terme $u$ de $t$ est dit {\it en
position d'application} ssi il existe un terme $v$ tel que $(u)v$ est
un sous-terme de $t$.

2) Un sous-terme $u$ de $t$ est dit {\it en position d'argument} s'il
n'est pas en position d'application.

\begin{lemme}
Soient $t$ un $\l$-terme $\b$-normal contenant $\a$. Si $\a : O , \G
\v_{\cal F} t : A$, alors $\a$ est en position d'argument dans $t$ et
donc, pour tout $\l$-terme $\b$-normal $u$, $t[u/ \a]$ est
$\b$-normal.
\end{lemme}

{\bf Preuve} Facile. \cqfd \\

{\bf D\'efinition} Un $\l$-terme $t$ est dit {\it bon} ssi le
$\l_{{\cal U V}}$-terme $t[{\cal U}_{X,X} /{\cal U}, {\cal V}_{X,X}
/{\cal V}]$ est bon.

\begin{lemme}
Si $t$ est un $\l$-terme bon $\b$-normal contenant l'une des
constantes ${\cal U}$ ou
${\cal V}$, alors $t[U /{\cal U}, V /{\cal V}]$ est normalisable et sa
forme normale contient $\a$.
\end{lemme}

{\bf Preuve} Si $t$ contient ${\cal U}$, alors ${\cal U}$ est
appliqu\'e \`a un seul terme $a$, et donc dans $t[U /{\cal U}, V
/{\cal V}]$ le sous-terme $(\l x \l d \langle x , \a \rangle)a$ se r\'eduit \`a $\l
d \langle a , \a \rangle$ sans cr\'eer des nouveaux redex. Si $t$ contient ${\cal
V}$, alors ${\cal V}$ est appliqu\'e \`a un terme de la forme
$(x)a_1...a_n$, et donc dans $t[U /{\cal U}, V /{\cal V}]$ le
sous-terme $(\l x (x) \a {\bf 1})(x) a_1...a_n$ se r\'eduit \`a $(x)
a_1...a_n \a {\bf 1}$ sans cr\'eer des nouveaux redex.  Dans les
deux cas la forme normale de $t[U /{\cal U}, V /{\cal V}]$ contient
$\a$. \cqfd

\begin{lemme}
Soient $u_1,...,u_n,v_1,...,v_m$ des $\l I$-termes, $A$ un type
contenant $X$ comme variable libre. Si $T = (({\cal
J}_{A,X})(x)u_1...u_n)v_1...v_m$ est typable, alors la forme normale
de $T$ contient $\a$.
\end{lemme}

{\bf Preuve} Consid\'erons le $\l$-terme $T' = (({\cal
J}'_{A,X})(x)u_1...u_n)v_1...v_m$.  D'apr\`es le lemme 3.6, $T'$ est
un $\l I$-terme. Comme $T'$ est fortement normalisable (sinon $T$ ne
le sera pas) et $A$ contient $X$, alors, d'apr\`es le lemme 3.6 et le
th\'eor\`eme 2.1, la forme normale $T'_1$ contient ${\cal U}$ ou
${\cal V}$. Soit $T'' = (({\cal
V}_{A,X})(x)u_1...u_n)v_1...v_m$. D'apr\`es le lemme 3.9, $T''$ est
fortement normalisable, et si on note $T''_1$ sa forme normale, alors
$T'_1 = \c{T''_{1}}$. $T''$ est un $\l_{{\cal U V}}$-terme bon, alors,
d'apr\`es le corollaire 3.1, $T''_1$ est bon ne contenant que des
${\cal U}_{X,X}$ et des ${\cal V}_{X,X}$ , donc $T'_1$ est bon. Donc,
d'apr\`es le lemme 3.13, $T'_1[U /{\cal U}, V /{\cal V}]$ contient
$\a$. Or comme $T = T'[U /{\cal U}, V /{\cal V}],$ alors la forme
normale de $T$ contient aussi $\a$. \cqfd \\

{\bf D\'efinition} On distingue trois r\`egles d'\'elimination
du quantifiacteur $\q$ :
\begin{center}
$(\q_e)_1 \; \F{ \G\v_{\cal F} \l x t :\q XA} { \G\v_{\cal F} \l x t: A[G/X] }$  \\[0,5cm] 

$(\q_e)_2 \; \F{ \G\v_{\cal F} (x) t_1 ...t_n :\q XA} { \G\v_{\cal F} (x) t_1 ...t_n
: A[G/X] }$ \\[0,5cm] 

$(\q_e)_3 \; \F{ \G\v_{\cal F} (\l x u)vt_1 ...t_n :\q XA} { \G\v_{\cal F} (\l x u)vt_1 ...t_n
: A[G/X] }$ 
\end{center}

\begin{lemme}
Soit $t$ un $\l I$-terme $\b$-normal. Si dans un typage $\G \v_{\cal
F} t : C$ on utilise la r\`egle $(\q_e)_2$, alors il existe un
$\l$-terme $\b$-normal $t'$ contenant $\a$ et $\G , \a : O \v_{\cal F}
t' : C$. \end{lemme}

{\bf Preuve} Dans le typage de $t$ on a donc 
\begin{center}
$\F{\F{\F{
\D \v_{\cal F} (x)u_1...u_n : \q X A}
{\D \v_{\cal F} (x)u_1...u_n : A[G/X]}}
{\matrix{.\cr.\cr.\cr}}}
{\D' \v_{\cal F} (x)u_1...u_nv_1...v_m : B}$
\end{center}
et dans $t$,
   $(x)u_1...u_nv_1...v_m$ n'est pas appliqu\'e \`a un autre $\l$-terme.\\

Consid\'erons le typage suivant 
\begin{center}
$\F{\F{\F{\F{
\D \v_{\cal F} (x)u_1...u_n : \q X A }
{\D \v (x)u_1...u_n : A[G^{\circ}/X]}
\quad
\F{}{\a : O \v_{\cal F} {\cal J}_{A,X} : A[G^{\circ}/X] \f A[G/X]}}
{\D , \a : O \v ({\cal J}_{A,X})(x)u_1...u_n : A[G/X]}}
{\matrix{.\cr.\cr.\cr}}}
{\D', \a : O \v (({\cal J}_{A,X})(x)u_1...u_n)v_1...v_m : B}$
\end{center}
Si $u$ est la forme normale de $(({\cal J}_{A,X})(x)u_1...u_n)v_1...v_m$, 
alors $\D', \a : O \v_{\cal F} u : B$ et, d'apr\`es le lemme 3.14, $u$
contient $\a$. Si $t = C \langle (x)u_1...u_nv_1...v_m \rangle$, alors posons $t' = C \langle u \rangle$. On a 
$t'$ est un $\l$-terme $\b$-normal contenant $\a$ et, d'apr\`es le
lemme 2.3,  $\G , \a : O \v_{\cal F} t' : C$.\cqfd

\begin{lemme}
Soient $D$ un $I$-type et $t$ un $\l I$-terme $\b$-normal. Si
$\v_{\cal F} t : D$, alors la r\`egle $(\q_e)_2$ ne peut pas \^etre
utilis\'ee.
\end{lemme}

{\bf Preuve} Sinon, d'apr\`es le lemme 3.15, il existe un $\l I$-terme
$\b$-normal $t'$ contenant $\a$ tel que $\a : O \v_{\cal F} t' :
D$. D'apr\`es le lemme 3.12, $\a$ est en position d'agument dans le
$\l$-terme $t'$, et donc $t''=t'[{\bf 1}/ \a]$ est un $\l K$-terme
$\b$-normal tel que $\v_{\cal F} t'' : D$. Contradiction.\cqfd \\

{\bf D\'efinition} On d\'efinit deux ensembles de type du syst\`eme ${\cal F}$ : $\q^+$ (ensemble
de types {\it $\q$-positifs}), et $\q^-$ (ensemble de types
{\it$\q$-n\'egatifs}) de la mani\`ere suivante :
 
 - Si $A = X$, alors $A \in \q^+$, et $A \in \q^-$ ;

 - Si $T \in \q^+$, et $T' \in \q^-$, alors, $T' \f T \in \q^+$, et $T
   \f T' \in \q^-$ ; 

 - Si $T \in \q^+$, et $X$ est libre dans $T$, alors $\q X T \in \q^+$.

Donc, $T$ est un type $\q^+$ (resp. $\q^-$) ssi le quantificateur du second ordre est
positif (resp. n\'egatif) dans $T$. 

\begin{lemme}
Soit $t$ un $\l I$-terme $\b$-normal tel que
$Fv(t)=\{x_1,...,x_n\}$. Si $t$ est terme $\eta$-long de type $A$ dans
le constext $\G = x_1 : A_1, ..., x_n : A_n$ sans utiliser la r\`egle
$(\q_e)_2$, alors les $A_i \in \q^-$ ($1 \leq i \leq n$) et $A \in
\q^+$.
\end{lemme}

{\bf Preuve} Par induction sur $t$.

- Si $t$ est une variable, alors $n = 1$, $A = X$ et $x_1 : A_1
\v_{\cal F} x_1 : X$. Donc $A_1 = X$ et par cons\'equent $A$ est
$\q^+$ et $A_1$ est $\q^-$.

- Si $t = (x_i)t_1...t_m$ ($m \geq 1$), alors $A = X$, $A_i =
B_1,B_2,...,B_m \f X$, et $\G \v_{\cal F} t_j : B_j$. Comme $t$ est un
$\l I$-terme alors, en appliquant l'hypoth\`ese d'induction sur les
$t_j$, on peut d\'eduire que les $A_k$ ($k \not = i$) sont $\q^-$ et
les $B_j$ sont $\q^+$. Donc $A_i$ est aussi $\q^-$.

- Si $t = \l x u$, alors $A = \q \sur{X} (B \f C)$ et $\G , x : B \v_{\cal F} u : C$. Par
hypoth\`ese d'induction, on d\'eduit que les $A_i$ et $B$ sont $\q^-$ et $C$ est $\q^+$. Donc $A$
est $\q^+$. \cqfd\\

{\bf D\'efinitions} 1) Si $D$ est un type du syst\`eme ${\cal F}$ alors on
note $\L(D) = \{ t$ $\b\eta$-normal clos /  $\v _{\cal F} t : D \}$.

2) Un type $D$ du syst\`eme ${\cal F}$ est dit {\it d\'emontrable} ssi
$card(\L(D)) \geq 1$.

\begin{theoreme}
Si $D$ est un $I$-type d\'emontrable du syst\`eme ${\cal F}$, alors $D
\in \q^+$.
\end{theoreme}

{\bf Preuve} Soit $t$ un $\l I$-terme $\eta$-longue de type
$D$. D'apr\`es le lemme 3.16, dans le typage de $t$, on n'a pas
utilis\'e la r\`egle $(\q_e)_2$. Donc, d'apr\`es le lemme 3.17, $D \in
\q^+$. \cqfd

\section{R\'esultats suppl\'ementaires}

\subsection{Les types entr\'ees et les types sorties}

 Commencant tout d'abord par d\'efinir les syst\`emes de typage ${\cal
F}_0$ et ${\cal F}_1$.\\

{\bf D\'efinitions} 1) Le syst\`eme de typage ${\cal F}_0$ est obtenue
en \'eliminant la r\`egle de typage $(\q_e)$ du syst\`eme ${\cal F}$.

2) Le syst\`eme de typage ${\cal F}_1$ est obtenue en remplacant la
r\`egle de typage $(\q_e)$ du syst\`eme ${\cal F}$ par
\begin{center}
$(\q_e)^1 \; \F{ \G\v_{{\cal F}_1} t :\q XA \;\;\; Y {\rm \; est \;
une \; variable}} { \G\v_{{\cal F}_1} t: A[Y/X] }$  
\end{center}

Les auteurs ont donn\'e des d\'efinitions (par typage) des types
entr\'ees et des types sorties du syst\`eme ${\cal F}$ (voir [2]).\\

{\bf D\'efinitions} 1) Un type clos $E$ du syst\`eme ${\cal F}$ est dit
{\it type entr\'ee} ssi si $\v_{\cal F} t : E$ alors $\v_{{\cal F}_0}
t : E$. Intuitivement un type entr\'ee est un type dont toutes les
d\'emonstrations se font dans le syst\`eme ${\cal F}_0$, et donc le
probl\`eme de typage pour ce type est d\'ecidable.

2) Un type clos $S$ du syst\`eme ${\cal F}$ (ne contenant pas une
constante de type fix\'es $O$) est dit {\it type sortie} ssi si $t$
est un $\l$-terme $\b$-normal tel que $\a : O \v_{\cal F} t : S$ alors
$\a \not \in Fv(t)$. Intuitivement cela veut dire que les fonctions
\`a valeurs dans un type sortie ind\'ependamment du type de leurs
arguments sont les fonctions constantes.\\

On a les r\'esultats suivants (voir [2] et [5]) :

\begin{theoreme}
1) Un type clos et $\q^+$ est un type entr\'ee et 
   sortie.

2) Un type entr\'ee est un type sortie.
\end{theoreme}

On ne sait pas si la r\'eciproque du 2) du th\'eor\`eme pr\'ec\'edent est
vraie. Nous allons voir que le lemme 3.15 apporte une r\'eponse
partielle \`a cette question : \\

En effet si un type clos $D$ n'est pas un type entr\'ee alors il
existe un $\l$-terme $\b$-normal clos $t$ tel que $\v_{\cal F} t : D$
et $\not \v_{{\cal F}_0} t : D$ c.\`a.d la r\`egle $(\q_e)$ est
utilis\'ee dans le typage. Le lemme 3.15 montre que si c'est la
r\`egle $(\q_e)_2$ qui est utilis\'ee et si $t$ est un $\l I$-terme,
alors $D$ n'est pas un type sortie.

\subsection{Les types simples}

On va montrer qu'on peut limiter l'\'etude des $I$-types \`a une
classe de types tr\`es simples.\\

{\bf D\'efinitions} 1) Les {\it types} du syst\`eme {$\cal S$}
(appel\'e syst\`eme simple) sont construits \`a partir des variables
de type en utilisant uniquement le connecteur $\f$.

2) Les r\`egles de typage du syst\`eme ${\cal S}$ sont $(ax)$,
$(\f_i)$ et $(\f_e)$.

3) On \'ecrit $\G \v_{\cal S} t : A$ si $t$ est typable sans le syst\`eme
{$\cal S$} de type $A$ dans le contexte $\G$.\\

Le th\'eor\`eme 2.3 reste valable dans le syst\`eme $\cal S$.\\

{\bf D\'efinition}
Si $A$ est un type du syst\`eme ${\cal F}$, alors on note ${\cal
S}_A$ l'ensemble des types du syst\`eme ${\cal S}$ obtenus en effacant
dans $A$ tous les quantificateurs (on suppose que toutes les variables
du type utilis\'ees dans $A$ ont des noms diff\'erents).\\

On a d\'emontr\'e dans [5] le r\'esultat suivant :

\begin{theoreme}\label{F->S} 
Soient $A$ un type $\q^+$ du syst\`eme ${\cal F}$, $A^* \in {\cal
S}_A$, et $t$ un $\l$-terme $\beta$-normal. Si $\v_{\cal F} t : A$,
alors $\v_{\cal S} t : A^*$
\end{theoreme}

Donc pour v\'erifier si un type clos $D$ de l'ensemble $\q^+$ est un
$I$-type, il suffit d'\'etudier cette question pour une cl\^oture d'un
type simple de l'ensemble ${\cal S}_D$. Ceci est plus facile \`a
\'etudier vu que le probl\`eme de typage pour les types simples est
d\'ecidable.\\

La r\'eciproque du th\'eor\`eme \ref{F->S} est \'evidement fausse si
on permet l'utilisation des types non propres. En effet, soient $A =
\q X \; \{(X \f X) \f (X \f \q Z X)\}$, $A^* = (X \f X) \f (X \f
X)$. On a $\v_{\cal S} {\bf id} : A^*$ mais $\not \v_{\cal F} {\bf id}
: A$. Remarquons que l'exemple pr\'ec\'edent ne marche pas si on
consid\`ere des $\l$-termes $\eta$-longs de type $A^*$. On va
pr\'esenter un contre exemple \`a la r\'eciproque du th\'eor\`eme
\ref{F->S} en se limitant aux $\l$-termes $\eta$-longs.\\

{\bf D\'efinition} Soient :

$T = \q X \{[ \q Y  (((Y \f \q Z  ((X,Y \f Z) \f Z)) \f X) \f
X) \f X] ,  X \f X\}$

$T^* = [(((Y \f ((X,Y \f Z) \f Z)) \f X) \f X) \f X] ,  X \f X$

$t = \l x \l y (x) \l z (x) \l u (z) \l v \l w ((w)(u) \l d \langle
y,v \rangle)v$

\begin{theoreme}
Le $\l$-terme $t$ est $\eta$-long de type $T^*$ mais $\not \v_{\cal F} t : T$.
\end{theoreme}

{\bf Preuve} Supposons que $\v_{\cal F} t : T$, alors il existe un
typage de $t$ o\`u les variables li\'ees de $t$ sont d\'eclar\'ees de
la mani\`ere suivantes :

$x :  \q Y \; (((Y \f \q Z \; ((X,Y \f Z) \f Z)) \f X) \f
X) \f X$

$y : X$

$z : (Y \f \q Z \; ((X,Y \f Z) \f Z)) \f X$

$u : (Y' \f \q Z \; ((X,Y' \f Z) \f Z)) \f X$ ($Y' \neq Y$)

$v : Y$

$w : X , Y \f Z$

$d : Y'$

et le sous-terme  $\langle y,v \rangle$ est de type $\q Z \;
((X,Y' \f Z) \f Z))$. D'o\`u $v$ est typable de type
$Y'$. Contradiction.\cqfd

\subsection{Quelques exemples}

{\bf D\'efinition} Soit $n \in \N$. Un $I$-type $D$ est dit d'ordre $n$
(resp. infini) ssi Card($\L(D)$) = $n$ (resp. $\L(D)$ est infini).\\

On va donner des exemples de $I$-types d'ordre quelconque.

\begin{theoreme}
$Id$ est $I$-type d'ordre $1$.
\end{theoreme}

{\bf Preuve} Facile. \cqfd \\

{\bf Notations} 1) Soient $A_1,...A_n$ des types du syst\`eme ${\cal
F}$. On note $A_1 \et ... \et A_n$ le type $\q X \{ (A_1,...,A_n \f X)
\f X\}$ o\`u $X$ ne figure pas dans les types $A_1,...,A_n$.

2) Soient $t_1,...,t_n$ des $\l$-termes. On note $\langle t_1,...,t_n
   \rangle$ le $\l$-terme $\l x (x) t_1...t_n$.\\
 
{\bf D\'efinition} Soit $n \geq 2$. On note :

$B_n = \q X \q Y_1 ... \q Y_n \{(Y_1 \f X),...,
(Y_n \f X) \f $

$[Y_1 ,  ... , Y_n \f (X \et Y_1 \et...\et Y_n)] \et (Y_1 \f X) \et ... \et (Y_n \f X)\}$

et, pour tout $1 \leq i \leq n$, 
$T_i = \l x_1 ... \l x_n \langle\l y_1 ... \l y_n  \langle (x_i)y_i
, y_1 , ... ,  y_n\rangle , x_1 , ... ,x_n \rangle$.

\begin{theoreme}
$B_n$ est un $I$-type d'ordre $n$ et $\L(B_n) = \{T_1, ... , T_n \}$.
\end{theoreme}

{\bf Preuve} Facile. \cqfd \\ 

{\bf Notation} Soient $x,y,z$ des variables et $n$ un entier. On note
$\l u[(x)(y)]^n(z)u$ la forme $\b\eta$-normale du $\l$-terme $\l u
(x)(y)...(x)(y)(z)u$ ($(x)(y)$ est r\'ep\'et\'e $n$ fois).\\

{\bf D\'efinition} Soit $B_{\infty} = \q X \q Y \{(X \f Y) , (Y
\f X) \f  [(X \f Y) \et  (Y \f X)]\}$ et, pour tout $(i,j) \in \N^2$,
$T_{i,j} = \l x \l y \langle \l u [(x)(y)]^i (x)u,\l u [(y)(x)]^j(y)u  \rangle$.

\begin{theoreme}
$B_n$ est un $I$-type infini et $\L(B_{\infty}) = \{T_{i,j}$ ; $ (i,j) \in \N^2\}$.
\end{theoreme}

{\bf Preuve} Facile. \cqfd \\ 

D'apr\`es le th\'eor\`eme 4.5 le type $B_2$ est un $I$-type qui
repr\'esente les bool\'eens. Il contient cinq quantificateurs. Une
question se pose : ``Peut-on trouver des $I$-types plus simples (avec
moins de variables ou de quantificateurs) pour les bool\'eens ?''
Nous allons montrer que les $I$-types contenant au plus deux
quantificateurs sont au plus d'ordre 1.

 \begin{theoreme} Si $D$ est un $I$-type d\'emontrable du syst\`eme
${\cal F}$ contenant un seul quantificateur, alors $D = \q X \{ X \f X
\}$.
\end{theoreme}

{\bf Preuve} Soient $D = \q X \{ A_1,...,A_n \f X\}$ un tel $I$-type
et $t$ un terme $\eta$-long minimal tel que $\v_{\cal F} t : D$. Alors
$t = \l x_1 ... \l x_n u$, $u = (x_i)u_1...u_m$ $(m \geq 0)$, et $x_1
: A_1 , ... , x_n : A_n \v_{\cal F} u : X$. Si $m \not = 0$, alors
$x_1 : A_1 , ... , x_n : A_n \v_{\cal F} u_1 : B_1,...,B_k \f
X$. D'apr\`es le choix de $t$, on a $k \geq 1$ et $\l x_1 ... \l x_n
((x_i) \l y_1 ... \l y_k u) u_2 ... u_m$ qui n'est pas un $\l I$-terme
est de type $D$. D'o\`u $m = 0$, $u = x_i$ et donc $n = 1$, et $D =\q
X \{X\f X\}$. \cqfd \\

{\bf Notations} 1) Si $B,A$ sont des types, alors, pour tout $n \geq
1$, on note $B^n \f A$ le type $B,...,B \f A$ o\`u $B$ est
r\'ep\'et\'e $n$ fois. 

2) Si $X$ et $Y$ sont deux variables de types, alors on note les
formules $\q X \q Y A$ et $\q Y \q X A$ par $\q X,Y A$.

\begin{theoreme}
Si $D$ est un $I$-type d\'emontrable du syst\`eme ${\cal F}$ contenant
deux quantificateurs, alors $D =\q X, Y \{ [(Y \f Y)^n \f X] \f X\}$,
$D = \q X \{ (\q Y (Y \f Y) \f X) \f X \}$, $D =\q X, Y \{ Y,(Y^n \f
X) \f X\}$ ou $D = \q X, Y \{(Y^n \f X) , Y \f X\}$.
\end{theoreme}

{\bf Preuve} Soient $D = \q X \q Y \{A_1,...,A_n \f X\}$ un $I$-type
d\'emontrable du syst\`eme ${\cal F}$ et $t$ un terme $\eta$-long
minimal tel que $\v_{\cal F} t : D$. Alors $t = \l x_1 ... \l x_n u$,
$u = (x_i)u_1...u_k$ $(k \geq 1)$, et $x_1 : A_1 , ... , x_n : A_n
\v_{\cal F} u : X$. Les $u_i$ ($1 \leq i \leq k$) ne peuvent pas
\^etre de type $B_1,...,B_m \f X$. En effet si c'est le cas, on
contredit soit la minimalit\'e de $t$ soit le fait que $D$ est un
$I$-type. Cherchons donc les termes $\eta$-longs minimaux $v$ tels que
$x_1 : A_1 , ... , x_n : A_n \v_{\cal F} \l y_1 ... \l y_m v :
B_1,...,B_m \f Y$ c.\`a.d. $\G = x_1 : A_1 ,... , x_n : A_n , y_1:B_1,
..., y_m : B_m \v_{\cal F} v : Y$. Supposons que $v = (z)v_1...v_l$
$(l\geq 0)$. Remarquons d'abord que les $v_i$ ($1 \leq i \leq l$) ne
peuvent pas \^etre ni de type $C_1,...,C_r \f Y$ ni de type
$D_1,...,D_{r'} \f X$ avec $(r' \geq 1)$. Donc ils sont forcemant de
type $X$. Soit $w$ un terme $\eta$-long minimal tel que $\G \v_{\cal
F} w : X$. $w$ est \'evidement une variable. On a deux cas \`a voir :

 -- Si $w = x_j$, alors $A_j = X$ et $D = \q X \q Y \{ X \f X
\}$. Contradiction.

 -- Si $w = y_j$ alors $x_1 : A_1 , ... , x_n : A_n \v_{\cal F}
\l y_1 ...\l y_m v[u/y_j] : B_1,...,B_m \f Y$ et $D$ n'est pas un
$I$-type.

Donc $v$ est une variable. De nouveau on a deux cas \`a voir :

 -- Si $v = x_j$, alors $A_j = Y$. Dans ce cas tous les $u_i$
($1 \leq i \leq k$) sont \'egaux \`a $x_j$. Donc $t = \l x \l y (x)
y...y$ ou $t = \l x \l y (y) x...x$. D'o\`u $D = \q X \q Y \{ Y,(Y^n
\f X) \f X \}$ ou $D = \q X \q Y \{(Y^n \f X) , Y \f X \}$.

 -- Si $v = y_j$ alors $B_j = Y$. Dans ce cas tous les $u_i$
($1 \leq i \leq k$) sont \'egaux \`a $\l z z$. D'o\`u $t = \l x (x) \l
z z...\l z z$ et $D = \q X \q Y \{( ( Y \f Y)^n \f X) \f X \}$.

Si $D = \q X \{ A_1,...,A_i \f \q Y (A_{i+1},...,A_n \f X) \}$, alors,
en reprenant la preuve que nous venons de faire, on d\'eduit que $i=0$
et donc $D$ est l'un des types trouv\'es pr\'ec\'edement.

Si $D = \q X \{ A_1 ,..., A_k, ..., A_n \f X \}$ et $A_k$ contient le
deuxi\`eme quantificateur, alors soit $t$ est terme $\eta$-longue
minimal tel que $\v _{\cal F} t : D$. Alors $t = \l x_1 ... \l x_n
(x_i)u_1...u_m$ et $x_1 : A_1 , ... , x_n : A_n \v_{\cal F}
(x_i)u_1...u_m : X$.  Les $u_i$ ($1 \leq i \leq m$) ne peuvent pas
\^etre de type $B_1,...,B_r \f X$. Donc $i = k$, $m = 1$ et $A_k = \q
Y(C_1 ,..., C_l \f Y) \f X$. En reprenant la preuve que nous avons
faite pr\'ec\'edement, on d\'eduit que $l=1$ et $C_1 = Y$. D'o\`u $D =
\q X \{ (\q Y (Y \f Y) \f X) \f X \}$. \cqfd

\begin{corollaire}
Si $D$ est un $I$- type du syst\`eme ${\cal F}$ contenant au plus deux
quantificateurs, alors $Card(\L(D)) \leq 1$.
\end{corollaire}

{\bf Preuve} D'apr\`es les th\'eor\`emes 4.7 et 4.8. \cqfd

\end{document}